\documentclass[11pt]{amsart}
\usepackage[mathscr]{eucal}
\usepackage{graphicx}
\newcommand{\beqa}{\begin{eqnarray}}
\newcommand{\eeqa}{\end{eqnarray}}
\newcommand{\beqan}{\begin{eqnarray*}}
\newcommand{\eeqan}{\end{eqnarray*}}

\newcommand{\beq}{\begin{equation}}
\newcommand{\eeq}{\end{equation}}
\newcommand{\beqn}{\begin{equation*}}
\newcommand{\eeqn}{\end{equation*}}

\newcommand{\R}{\mathbb{R}}
\newcommand{\PP}{\mathscr{P}}
\newcommand{\F}{\mathscr{F}}

\newtheorem{theo}{Theorem}

\newtheorem{lem}[theo]{Lemma}

\newtheorem{coro}[theo]{Corollary}

 \newcommand{\bit}{\begin{itemize}}
\newcommand{\eit}{\end{itemize}}
\newcommand{\argmin}{\arg\!\min}
\newcommand{\argmax}{\arg\!\max}

\begin{document}

\title{On some algorithmic aspects of hypergraphic matroids}
\date{\today}

\author[M. Ba\"{\i}ou]{Mourad Ba\"{\i}ou}
\address[M. Ba\"{\i}ou]{ CNRS,  and Universit\'e Clermont II,
Campus des c\'ezeaux BP 125, 63173 Aubi\`ere cedex, France.}

\email[M. Ba\"{\i}ou]{baiou@isima.fr}

\author[F. Barahona]{Francisco Barahona}
\address[F. Barahona]{IBM T. J. Watson research Center, Yorktown Heights, NY 10589, USA. }
\email[F. Barahona]{barahon@us.ibm.com}

\begin{abstract}
Hypergraphics matroids were studied first by Lorea \cite{Lo} and later
by Frank et al \cite{Frank}. They can be seen as generalizations of graphic
matroids. Here we show that several algorithms developed for the graphic
case can be extended to hypergraphic matroids. We treat the following:  
the separation 
problem for the associated polytope, testing independence,
separation of
partition inequalities, computing the rank of a set, computing the strength,
computing the arboricity and network reinforcement.
\end{abstract}

\keywords{Hypergraphic Matroids, Combinatorial Optimization, Network Reinforcement}



\maketitle

\section{Introduction}
Hypergraphic matroids were introduced by Lorea \cite{Lo} and later studied by Frank et al \cite{Frank}.
They showed that the notion of circuit-matroid of graphs can be generalized to hypergraphs.
In \cite{Frank}
they generalized the notion of spanning trees 
to hypertrees and extended a theorem of Tutte \cite{Tu}
and Nash-Williams \cite{NW}, 
to give the maximum number of disjoint hypertrees contained in a hypergraph. 

Different algorithms associated with 
graphic matroids have been developed based on the graph
structure. Here we extend some of these algorithms to the hypergraphic case.
More precisely, we give an algorithm for separating a vector
from the associated matroid polytope.
For the graphic case this had been treated in \cite{PQ} and \cite{PW}.
 We also show that when applying the
greedy algorithm to find a maximum weighted
independent set \cite{Ed}, testing independence at each iteration can be done
by finding a minimum cut in an associated graph.
We also treat partition inequalities,
this leads to an algorithm
to compute the rank of a set. For the graphic case, partition inequalities
were treated in \cite{Cu} and \cite{BaDom}.
The maximum number of disjoint spanning trees in a graph has been
proposed in \cite{Gu} as a measure of the {\it strength} of a network.
An algorithm for computing the strength
of a graph was given in \cite{Cu}, and later in \cite{Gu2} and \cite{CC}.
Here we give an algorithm to compute a similar measure in a hypergraph, namely the
maximum number of disjoint hypertrees.
The {\it arboricity} of a graph is the minimum number of edge-disjoint forests
into which the edge-set can be decomposed. Algorithms to compute the 
arboricity were given in \cite{PQ},  \cite{GGT}, and \cite{Ga}. Here we give an algorithm
to compute a similar measure in hypergraphs.
Also the {\it reinforcement problem} was studied in \cite{Cu}, this consists of
given an initial graph $G$, a target value $k$ and a set of candidates
edges, find a minimum cost set of candidate edges to be added to $G$
so that the resulting graph has $k$ disjoint spanning trees.
Here we give an algorithm for a similar reinforcement problem in a
hypergraph.  

The problems mentioned above
can be solved with general matroid algorithms or submodular flows, using the 
characterizations in \cite{Frank} and \cite{Frank2}.
In particular, the reinforcement problem is equivalent to finding a minimum cost
sub-hypergraph that has a rooted $k$-edge-connected orientation, solvable 
via submodular flows, see \cite{Frank2}.
However the algorithmic details are not given in these references.
Here we give algorithms that reduce all these to a sequence of minimum cut problems
in an auxiliary graph.

This paper is organized as follows. In Section~\ref{prel} we give some 
definitions and notation. 
In Section~\ref{basic} we study the separation
problem and the greedy algorithm.
Section~\ref{separDom} is devoted
to partition inequalities and to computing the rank of a set.
The strength of a hypergraph is studied in Section~\ref{strength}.
The arboricity is studied in Section~\ref{arbo}.
In Section~\ref{reinf} we study the reinforcement problem.

\section{Preliminaries}\label{prel}
Let $H=(V, E)$ be a hypergraph. 
For a non-empty set $X \subset V$ and $F \subseteq E$,  $F[X]$ denotes the set of
hyperedges in $F$ contained in $X$.
Let $\PP=\{V_1, \ldots, V_k\}$ be a family of non-empty subsets of $V$ with 
$V_i\cap V_j=\emptyset$ for $i \ne j$, 
we denote by $\delta_F(\PP)$ the set of hyperedges in $F$ included in $\cup_i V_i$
and that intersect at least two sets in $\PP$. 
For a hypergraph $H'=(V,E')$ sometimes we use $\delta_{H'}(\PP)$ instead of
$\delta_{E'}(\PP)$. Also when there is no confusion we use $\delta(\PP)$ instead
of $\delta_F(\PP)$.

A linear system $Ax\leq b$, where $A$ and $b$ are rational, is called 
{\it totally dual integral (TDI)} if for any $c\in {\mathbb {Z}}^{n}$ such that there is an
optimal solution to the linear program
$\max \{c^{\mathrm {T} }x \,  \vert \, Ax\leq b \}$,
there is an integer optimal dual solution. 
Edmonds and Giles \cite{EG} showed that if a polyhedron $P$ is the solution set of a TDI system  $Ax\le b$, where $b$ has all integer entries, then every vertex of  $P$
 is integer-valued. 

For  a {\it multiset} $F \subseteq E$, its {\it incidence vector } $x^F \in \R^E$
is defined
as follows,
$x^F(e)$ is the {\it multiplicity} of $e$ in  $F$.
For a vector $x \in \R^E$, and
$S \subseteq E$ we use $x(S)$ to denote $\sum_{e \in S} x(e)$.

Let $D=(V,A)$ be a directed graph, for $S \subseteq V$, we denote by
$\delta^+(S)$ the set 
$\delta^+(S)=\{(u,v) \in A \,  \vert \, u \in S, v \notin S\}$.
Given
two distinguished  vertices
$s$ and $t$, for a set $S\subset V$, with $s\in S$,
$t \notin S$, the set of arcs
$\delta^+(S)$ is called an {\it $st$-cut}.
Given a capacity vector $c \in \R^A_+$, a {\it minimum} $st$-cut is an $st$-cut
$\delta^+(S)$ such that $c(\delta^+(S))$ is minimum.
A minimum $st$-cut
can be found in $O(|V|^3)$ time, see \cite{GT}.

A {\it hyperforest} in $H=(V,E)$ is a set $F \subseteq E$ such that
$|F[X]| \le |X|-1$ for every non-empty $X \subseteq V$. A hyperforest $F$ is called a 
{\it hypertree} of $H$ if $|F|=|V|-1$. 
If $H$ is a graph, 
$F \subseteq E$ is a hypertree if and only if $F$ is a spanning tree. It was proven by
Lorea \cite{Lo} that the hyperforests of a hypergraph form the family of independent sets
of a matroid. These are called {\it hypergraphic matroids}. 
Frank et al. \cite{Frank} further studied hypergraphic matroids.
Hypertrees are bases of
hypergraphic matroids.

\section{Basic algorithms}\label{basic}

Here we study the separation problem for polytopes of hypergraphic matroids, and 
the application of the greedy algorithm to find a maximum weighted independent
set. The algorithms in this section are quite simple, we present them for the sake
of completeness.

\subsection{Separating from the associated polytope}\label{sepPol}
For a hypergraph $H=(V,E)$,
let $P(H)$ be the convex hull of incidence vectors of independent sets
of the associated matroid.
Let $r(S)$ be the rank of $S$, for $S \subseteq E$.
A fundamental theorem of Edmonds \cite{Ed} gives an explicit description
of $P(H)$ as below.

\begin{theo}
$
P(H)=\{ x \in \R^E \,   \vert  \, x \ge 0, \ x(S) \le r(S), \, \forall S \subseteq E \}.
$
\end{theo}

The {\it separation problem} for $P(H)$ consists of 
given a vector $\bar x$,  
find a hyperplane separating $\bar x$ from $P(H)$, or
decide that $\bar x \in P(H)$. For the graphic case this has been
solved in \cite{PQ} and \cite{PW}.
A polynomial algorithm for the separation problem enables us to optimize
a linear function over $P(H)$ in polynomial time, cf. \cite{GLS}.
Frank et al \cite{Frank}
proved the following.

\begin{theo}\label{rank}
For $S \subseteq E$,
the rank of $S$, $r(S)$ is given by the following formula:
\beqn
r(S) = \min \{\vert V \vert- \vert \PP \vert + \vert \delta_S(\PP)\vert : \, \PP \mbox{ a partition of } V\}.
\eeqn
\end{theo}

We assume that $0 \le \bar x(e) \le 1$ for all $e \in E$, otherwise we have a 
separating hyperplane immediately.
Thus for $S \subseteq E$,
let us suppose that there is a partition $\PP=\{V_1, \ldots, V_k\}$ of $V$
such that $$\bar x(S) > |V| - |\PP | + |\delta_S(\PP)|.$$
Since $\bar x(e) \le 1$, for all $e \in E$, we have $\bar x(\delta_S(\PP)) \le |\delta_S(\PP)|$.
And $\bar x(S)=\bar x(\delta_S(\PP)) + \sum_i \bar x(S[V_i])$, therefore
$\sum_i \bar x(S[V_i]) > |V|-|\PP|=\sum_{i=1}^k (|V_i|-1)$. Thus $\bar x(S[V_r]) > |V_r|-1$
for at least one index $r$. This suggests to look for a 
non-empty node-set $W \subset V$
that maximizes $\bar x(E[W])-|W| + 1$. This can be reduced to a sequence of minimum cut problems
as follows. This construction is inspired on a construction of Lawler \cite{La}.

We define a directed graph $D=(V', A)$, with $V'=\{s, t\} \cup V \cup \{e', \, e'' \, : \, e \in E\}$.
The arc set is defined as follows:
\bit
\item[(a)]
For every node $v \in V$ we set an arc $(s,v)\in A$ with capacity equal to one. 
\item[(b)]  For every hyperedge $e \in E$ we define an arc $(e', e'')$ with capacity $\bar x(e)/2$,
also for each node $u \in e$ we define the arcs $(u,e')$ and $(e'', u)$ with infinite capacity.
We also define arcs $(e',t)$ and $(e'',t)$ with capacity $\bar x(e)/2$.
\item[(c)]  Finally we pick a node $\bar v \in V$ and add an arc $(\bar v, t)$ with infinite capacity. This is to avoid having $\delta^+(V' \setminus \{t\})$ as a minimum
$st$-cut.
\eit
We find a minimum $st$-cut, and we repeat for each $\bar v \in V$. 

Assume that $\{s\} \cup T$ defines one of these minimum cuts. Let $T'=V\cap T$,
$\bar T'=V \setminus T'$,
\hbox{$E_1=\{e \, | \, e' \in T, \, e'' \notin T\}$}, $E_2=\{e \, | \, e', e'' \in T\}$, 
$E_3=\{e\, | \, e', e'' \notin T\}$. The capacity of this cut is
$$|\bar T'| + \bar x(E_1) + \bar x(E_2)= |\bar T'| - \bar x(E_3) + \bar x(E).$$
Notice that $E_3$ is the set of hyperedges contained in $\bar T'$.
Thus for each $\bar v \in V$ we have minimized $|W| - x(E[W])$, with the constraint
that $\bar v \in W$. If the value of the minimum is less than one, we have
a set $W \subseteq V$ with $\bar x(E[W]) > |W|-1$.  Then with $S=E[W]$, 
$V \setminus W=\{v_1, \ldots, v_r\}$ and 
$\PP=\{W, \{v_1\},\ldots,\{v_r\}\}$, the inequality $\bar x(E[W]) > |W|-1$ is
equivalent to $\bar x(S) > |V| - |\PP | + |\delta_S(\PP)|$.  Then we have a separating
hyperplane.
Now we can state the following.

\begin{theo}
For hypergraphic matroids the separation problem for the associated  polytope
reduces to $|V|$ minimum cut problems in a graph with $O(|V| + |E|)$
nodes.
\end{theo}

\subsection{Finding a maximum weighted hyperforest}
For a weight vector $w \in \R_+^E$, we deal with the problem 
of finding an independent set $I$ such
that $w(I)$ is maximum.
Edmonds \cite{Ed} showed that this can be done
with the greedy algorithm below.

\bit
\setlength{\itemindent}{1.2cm}
\item[\bf Step 0.] Order the elements so that $w(e_1) \ge w(e_2) \ge \ldots w(e_m)$.
Set $I=\emptyset$.
\item[\bf Step 1.] For $i=1,\ldots,m$ do
\bit
\item[] if $I \cup \{e_i\}$ is independent, set $I \leftarrow I \cup \{e_i\}$.
\eit
\eit

Here we need to show how to decide whether the set $I'=I \cup \{e_i\}$ in Step 1,
remains independent.
This is equivalent to look for the minimum of
$$|\cup\{e \, | \, e \in F \}|-|F|$$ for $F \subseteq I'$,
with $e_i \in F$.
Then $I'$ is independent if and only if the minimum is
at least $1$. 

Now we have to see that this reduces to a minimum cut problem.
For that we build a directed graph with node-set $\{s,t\} \cup \{e \, | \, e\in I'\} \cup \{v \, | \, v \in V'\}$. Here $V'=\cup \{ e \, | \, e \in I' \}$.
The arc-set is as follows:
\bit
\item There is an arc $(s,e)$ for each $e \in I$, with capacity equal to one.
\item There is an arc $(v,t)$ for each $v \in V'$, with capacity equal to one.
\item There is an arc $(e,v)$ for each $e \in I'$, and for each $v \in e$, with infinite capacity.
\item There is an arc $(s, e_i)$ with infinite capacity.
\eit

Assume that $\{s\} \cup T$ defines  a minimum $st$-cut. Then  $F=T\cap I'$,
and $\cup\{e \, | \, e \in F \}=T \cap V'$.
If $C$ is the capacity
of this cut, then $C - |I'|=|\cup\{e \, | \, e \in F \}|-|F|$ is the 
minimum value needed. We can state the following.

\begin{theo}
Finding a maximum weighted hyperforest reduces to $|E|$ minimum
cut problems in a graph with $O(|V|+|E|)$ nodes.
\end{theo}

\section{Partition Inequalities and Computing the Rank}\label{separDom}
Here we study the separation problem for partition inequalities.
This will be used to compute the rank of a set, and also it
will be used in the next two sections.
Frank et al \cite{Frank} proved the following.

\begin{theo}\label{teo}
A hypergraph contains $k$ disjoint hypertrees if and only if 
$$
|\delta(\PP)| \ge k (|\PP| - 1)
$$
holds for every partition $\PP$ of $V$.
\end{theo}

Thus for a hypergraph $H=(V,E)$, we can write the system of inequalities below,
that must be satisfied by the incidence vectors of hypertrees.
\beqa
&&x(\delta(\PP)) \ge |\PP| -1, \quad
\mbox{for all partitions $\PP$ of $V$,}\label{in1}\\
&& x \ge 0.\label{in2}
\eeqa

Let $P$ be the polyhedron defined by \eqref{in1}--\eqref{in2}.
It can have extreme points with components
greater than one. For instance if $H=(V=\{v_1, v_2, v_3\}, E=\{e_1, e_2\})$ and
$e_1=\{v_1, v_2, v_3\}$, $e_2=\{v_1, v_2, v_3\}$, then $P=\{(x_1, x_2) \, : \,
x_1 + x_2 \ge 2, \, x \ge 0\}$. The extreme points are $(2,0)$ and $(0,2)$.
If $H$ is a graph, this is an unbounded
 polyhedron whose extreme points are the incidence
vectors of spanning trees. This follows from a theorem of Tutte \cite{Tu} and
Nash-Williams \cite{NW}. Theorem~\ref{teo} extends this to hypergraphs.

\begin{coro}
The polyhedron $P$ has integral extreme points.
\end{coro}
\proof
First notice that the incidence vector of any hypertree satisfies \eqref{in1}-\eqref{in2}.
Let $\bar x$ be an extreme point of $P$.
This is a rational vector, so there is a positive integer $k$ such that $k \bar x$ is an integer
vector. Consider the hypergraph $\bar H=(V, \bar E)$, where for each $e \in E$ there are
$k \bar x(e)$ copies of $e$ in $\bar E$. Theorem~\ref{teo} implies that $\bar H$ contains
$k$ disjoint hypertrees.
Thus there are $k$ hypertrees $T_1,\ldots, T_k$ such that
$k\bar x \ge (x^{ T_1} + \ldots + x^{ T_k})$. Here each edge $e$
is being used at most $k \bar x(e)$ times, and it could be used
more than once in some hypertrees.
Thus $\bar x \ge \frac{1}{k}(x^{ T_1} + \ldots + x^{ T_k})$.

Since $\bar x$ is an extreme point, we have 
$\bar x = \frac{1}{k}(x^{ T_1} + \ldots + x^{ T_k})$, and
 $x^{ T_1} = \ldots = x^{ T_k}$.
\endproof

Now we discuss the separation problem for inequalities \eqref{in1}. 
We assume that 
$\bar x \in \R^E_+$ is an input vector. We are going to solve
\beq \label{fundSep}
\mbox{minimize } \quad \bar x(\delta(\PP)) - \beta (|\PP | - 1),
\eeq
where the minimum is taken among all partitions $\PP$ of $V$. Since
the partition $\PP=\{V\}$ is also possible, the minimum is always less
than or equal to zero. In this section we need $\beta=1$, but in the following
sections this number could be different.

We fix an arbitrary node $r \in V$, and for each non-empty $S \subseteq V$ we
define
\beq \label{func}
 f(S)=
 \begin{cases}
 \beta +\bar x(E[S]), \mbox{ if $r \notin S$,}\\
 \bar x(E[S]), \mbox{ if $r \in S$.}\\
 \end{cases}
 \eeq
 
The function $f$ is intersecting supermodular, i.e.,
$f(S \cup T) + f(S \cap T) \ge f(S) + f(T)$ for $S, T \subseteq V$, with
$S \cap T \ne \emptyset$.
We associate a variable $y(v)$ to each node $v \in V$,
and 
we propose to solve the
linear program
  \beqa
 &&\min y(V)\label{lp1}\\
 &&y(S) \ge f(S) \mbox{ for $\emptyset \ne S \subseteq V$.} \label{lp2}
 \eeqa
 Edmonds \cite{Ed2} showed that this
 can be solved with the greedy algorithm below.
 \bit
 \setlength{\itemindent}{1.2cm}
 \item[\bf Step 0.] Set $\bar y(v)=\beta +\bar x(E)$, for all $v \in V$; $\F = \emptyset$.
 \item[\bf Step 1.] Pick a node $\bar v$ that does not belong to a set in $\F$.
 If such a node does not exist stop, otherwise let 
 $\bar S \in \argmin \{ \bar y(S) - f(S)\, | \, \bar v \in S \}$,
 and let $\alpha = \bar y( \bar S)-f(\bar S)$.
 \newline
 Set $\bar y(\bar v) \leftarrow \bar y(\bar v) - \alpha$.
 $\F \leftarrow \F \cup \{\bar S\}$.
 \item[\bf Step 2.] Uncross: while there are two sets $S$ and $T$ in $\F$ with
 $S \cap T \ne \emptyset$, do \linebreak
 $\F \leftarrow \F \setminus \{S, T\} \cup \{ S \cup T \}$.
 Go to Step 1.
 \eit
 The vector $\bar y$ is built so it satisfies \eqref{lp2}. 
 A set $S \subseteq V$ is called {\it tight} if $\bar y(S)=f(S)$.
 It follows from the supermodularity of $f$ that if $S$ and $T$ are tight
 and $S \cap T \ne \emptyset$, then $S \cup T$ is also tight.
 This justifies the uncrossing operation in Step 2.
 At the end the family $\F$ defines a 
 partition of $V$ so that $\bar y(S)=f(S)$ for every $S \in \F$. 
 Now we discuss the two possible values of the solution.
 \bit
 \item
 If the value of the
 optimum is $\bar x(E)$ 
 we can take any partition $\{V_1, \ldots, V_p\}$ of $V$, add the 
 associated inequalities  \eqref{lp2}, and we have
 $$\bar x(E)=\bar y(V)=\sum \bar y(V_i) \ge \sum f(V_i)=\sum_{i=1}^p \bar x(E(V_i)) +
 \beta(p -1).$$
 Since $\bar x(E) = \sum \bar x(E(V_i)) + \bar x(\delta(V_1, \ldots, V_p))$, we obtain
 $$\bar x(\delta(V_1, \ldots, V_p)) \ge \beta(p-1).$$ This shows that inequalities \eqref{in1} are
 satisfied. 
 
 \item
 Now assume that the value of the optimum is greater than $\bar x(E)$. We set 
 $\bar z_S=1$ if $S \in \F$ and $\bar z_S=0$ otherwise, for each $S \subseteq V$. This
 vector is a feasible solution of the dual of \eqref{lp1}-\eqref{lp2}. We have
 $\bar y(V)=\sum\{\bar y(S) \, | \, S \in \F\}=\sum\{f(S) \, | \, S \in \F \}=
 \sum\{ f(S) \bar z_S \, | \, S \subseteq V \}$. This shows that $\bar z$ is an optimal
 dual solution. Let $\F=\{W_1, \ldots, W_p\}$.
 We have 
 \beqn
 \sum\{f(S) \, | \, S \in \F \}=
 \sum_i \bar x(E(W_i)) + \beta(p-1) > \bar x(E).
 \eeqn
 Then
 \beqn
 0 > 
 \bar x(E) - \sum_i \bar x(E(W_i)) -\beta(p-1)=\bar x(\delta(W_1,\ldots, W_p)) - \beta(p-1).
 \eeqn
 Since $\bar z$ is an optimal solution, we have 
 a most violated partition inequality.
 \eit
 
It remains to show how to find the set $\bar S$ in Step 1.  
For that
we define a directed graph $D=(V', A)$, with $V'=\{s, t\} \cup V \cup \{e', \, e'' \, : \, e \in E\}$. Let $\eta(v)=\bar y(v)$ for $v \ne r$, and $\eta(r)=\bar y(r) + \beta$. 
Let $\eta^+(v)=\max\{0, \eta(v) \}$, and $\eta^-(v)=\min \{0, \eta(v)\}$ for
$v \in V$.
The arc set is defined as follows.
\bit
\item
For every node $v \in V$ we set an arc $(s,v)\in A$ with capacity equal to $\eta^+(v)$,
and we set an arc $(v,t)\in A$ with capacity equal to $-\eta^-(v)$.
\item For every hyperedge $e \in E$ we define an arc $(e', e'')$ with capacity $\bar x(e)/2$,
also for each node $u \in e$ we define the arcs $(u,e')$ and $(e'', u)$ with infinite capacity.
We also define arcs $(e',t)$ and $(e'',t)$ with capacity $\bar x(e)/2$.
\item Finally we add an arc $(\bar v, t)$ with infinite capacity.
\eit

Assume that $\{s\} \cup S$ defines a minimum $st$-cut. Let $S'=V\cap S$,
$\bar S'=V \setminus S'$,
\hbox{$E_1=\{e \, | \, e' \in S, \, e'' \notin S\}$}, $E_2=\{e \, | \, e', e'' \in S\}$, 
$E_3=\{e\, | \, e', e'' \notin S\}$. The capacity of this cut is
\beqn
\eta^+( \bar S' ) -\eta^-( S')+ \bar x(E_1) + \bar x(E_2) = \eta(\bar S') - \bar x(E_3) + \bar x(E) -\eta^-(V).
\eeqn
Notice that $E_3$ is the set of hyperedges contained in $\bar S'$.
Thus if $\gamma$ is the capacity of the minimum cut obtained, then 
$\gamma -\beta - \bar x(E) + \eta^-(V)$ is the value $\alpha$ in
Step 1.
This minimizes $\bar y(W) - f(W)$, with the constraint
that $\bar v \in W$. Now we can state the following.

\begin{theo}
The separation problem for inequalities \eqref{in1} reduces to $|V|$ minimum
cut problems in a graph with $O(|V|+|E|)$ nodes.
\end{theo}

 Frank et al \cite{Frank} proved that for a set $F \subseteq E$ the rank
 of $F$ is given by
\beqn
\min \{ |\delta_F(\PP)| -|\PP | + |V| \, : \, \PP \mbox{ a partition of $V$}\}.
\eeqn
Thus to compute the rank of $F$ we just have to apply the separation algorithm above
 to the incidence vector of $F$. This leads to the following.
 \begin{coro}
 For hypergraphic matroids,
 computing the rank of a set $F \subseteq E$ reduces to $|V|$ minimum
cut problems,  in a graph with $O(|V|+|E|)$ nodes.
 \end{coro}

\section{Strength of a Network}\label{strength}

For a network represented by a graph $G$, the maximum number 
of edge-disjoint spanning trees contained in $G$, has been proposed as
a measure of the strength of the network, see \cite{Gu}. Algorithms to
compute the strength of a graph have been given in \cite{Cu}, \cite{Gu2} and \cite{CC}.
Here we give an algorithm to compute the maximum number of
disjoint hypertrees in a hypergraph.

Based on Theorem~\ref{teo},
for a hypergraph $H=(V,E)$ and a capacity vector $c \in \R^E_+$ we give
an algorithm to find the minimum of
\beqan
&&\frac{c(\delta(\PP))}{|\PP|-1},
\eeqan
among all partitions $\PP$ of $V$, with $|\PP| \ge 2$.
This gives the value of a maximum packing of hypertrees with
capacities $c(e)$ for each hyperedge $e \in E$. 
We follow an approach similar to the one of \cite{CC} for the graphic case.
We use 
Newton's method, cf. \cite{Rad}, as below.

\pagebreak[1]
\vskip 0.2cm
\centerline{{\bf Newton's method}}
\begin{itemize}
\setlength{\itemindent}{1.2cm}
\item[{\bf Step 0.}] 
Pick any partition $\bar \PP$ of $V$ with $|\bar \PP| \ge 2$.
Set \beqn
\beta=\frac{c(\delta(\bar \PP))}{|\bar \PP|-1}.
\eeqn
\item[{\bf Step 1.}] Find $\hat \PP=\argmin \Big\{
c(\delta(\PP)) - \beta(|\PP| - 1)
\Big\}$.
 \item[{\bf Step 2.}] If $c(\delta(\hat \PP)) - \beta(|\hat \PP| - 1)< 0$,
 update $\beta$ as
 \beqn
\beta=\frac{c(\delta(\hat \PP))}{|\hat \PP|-1}.
\eeqn
  and go to Step 1.
\newline
Otherwise  $c(\delta(\hat \PP)) - \beta(|\hat \PP| - 1)=0$,
and we stop.
\eit
 The minimum in Step 1 is found among all partitions $\PP$ of $V$ with $|\PP| \ge 2$.
 Notice that because of the definition of $\beta$, in Step 2 we always have
 $c(\delta(\hat \PP)) - \beta(|\hat \PP| - 1) \le 0$.
 If $\PP_1, \ldots, \PP_k$ is the sequence of partitions obtained we have 
 $c(\delta(\PP_i))/(|\hat \PP_i| - 1) > c(\delta(\PP_{i+1}))/(|\hat \PP_{i+1}| - 1) $,
 for $i=1, \ldots, k-2$. It follows from results of \cite{Rad}, (Section 3, Lemma 1),
  that $|\PP_i| > |\PP_{i+1}|$, 
 for $i=2,\ldots, k-2$. Thus this algorithm converges
 in at most $|V|$ iterations. 
 
Finding the minimum in Step 1 is similar  to the separation problem of
inequalities \eqref{in1} in Section~\ref{separDom}. The only difference
is the number $\beta$ that multiplies $|\PP|-1$. 
Notice that after each iteration of Newton's
method the number $\beta$ is decreasing. 
In the algorithm of Section~\ref{separDom},
for each node $v$ through the
iterations of Newton's method
we have to solve 
a sequence of minimum $st$-cut problems where for each node $w$
the capacities of 
the arcs $(s,w)$ are not increasing, and the capacities of the arcs $(w,t)$
are not decreasing, while the capacities of all the other arcs remain the same.
This suggests the use of the parametric
network flow algorithm of \cite{GGT}. 
There they deal with the case where for each node $w$ the capacities
of the arcs $(s,w)$ are not decreasing, and the capacities of the arcs
$(w,t)$ are not increasing. One can reverse the orientation of
every arc and look for a minimum $ts$-cut in the new graph.
Then each of these sequences can be
solved with the same asymptotic complexity as one application of
the push-preflow algorithm of \cite{GT}. This leads to the following.

\begin{theo}
The strength of a hypergraph can be computed with the same asymptotic complexity 
of $|V|$ applications of the push-preflow algorithm in a graph with
$O(|V|+|E|)$ nodes. 
\end{theo}

\section{Arboricity}\label{arbo}
 The {\it Arboricity} of a graph is the minimum number of edge-disjoint
 forests into which the edge-set can be decomposed. Nash-Williams \cite{NWarb}
 gave a characterization of this number,
Frank et al. \cite{Frank} extended this to hypergraphs as below.
 \begin{theo}
 A hypergraph $H=(V,E)$ can be partitioned into $k$ disjoint hyperforests if and only if
 for every $X \subseteq V$ 
 $$ |E[X]| \le k (|X| -1).$$
 \end{theo}
 
 Algorithms for computing the arboricity of a graph have been given in
 \cite{PQ}, \cite{Ga} and \cite{GGT}. Here we give an algorithm to compute the arboricity of
 a hypergraph. For that we need to compute
 \beqn
 \max \frac{|E[X]|}{|X|-1},
 \eeqn
for $X \subseteq V$ and $|X| \ge 2$. 
We assume that there is no hyperedge $e$ with $|e|=1$, otherwise $k=0$.
As in the last section we use Newton's
method as follows.

\pagebreak[4]
 \vskip 0.2cm
\centerline{{\bf Newton's method}}
\begin{itemize}
\setlength{\itemindent}{1.2cm}
\item[{\bf Step 0.}] 
Pick any  $\bar X \subseteq V$ with $|\bar X| \ge 2$.
Set \beqn
\beta=\frac{|E[\bar X]|}{|\bar X|-1}.
\eeqn
\item[{\bf Step 1.}] Find $\hat X =\argmax \Big\{
|E[X]| - \beta(|X| - 1)
\Big\}$.
 \item[{\bf Step 2.}] If $|E[\hat X]| - \beta(|\hat X| - 1)> 0$,
 update $\beta$ as
 \beqn
\beta=\frac{|E[\hat X]|}
{|\hat X|-1}.
\eeqn
  and go to Step 1.
\newline
Otherwise  $|E[\hat X]| - \beta(|\hat X| - 1)=0$,
and we stop.
\eit

The maximum in Step 1 is found among all $X \subseteq V$ with $|X| \ge 2$.
 Because of the definition of $\beta$, in Step 2 we always have
 $|E[\hat X]|- \beta(|\hat X| - 1) \ge 0$.
 If $X_1, \ldots, X_k$ is the sequence of node-sets obtained,
we have 
 $|E[X_i]|/(|X_i| - 1) > |E[X_i]|/(|X_{i+1}| - 1) $,
 for $i=1, \ldots, k-2$. If follows from the results of \cite{Rad}
 that
  $|X_i| > |X_{i+1}|$, 
 for $i=2,\ldots, k-2$. Thus this algorithm converges
 in at most $|V|$ iterations.

To find the maximum in Step 1, 
we define a directed graph $D=(V', A)$, with $V'=\{s, t\} \cup V \cup \{e', \, e'' \, : \, e \in E\}$.
The arc set is defined as follows:
\bit
\item[(a)]
For every node $v \in V$ we set an arc $(s,v)\in A$ with capacity equal to $\beta$. 
\item[(b)]  For every hyperedge $e \in E$ we define an arc $(e', e'')$ with capacity $1/2$,
also for each node $u \in e$ we define the arcs $(u,e')$ and $(e'', u)$ with infinite capacity.
We also define arcs $(e',t)$ and $(e'',t)$ with capacity $1/2$.
\item[(c)]  Finally we pick a node $\bar v \in V$ and add an arc $(\bar v, t)$ with infinite capacity. This is to avoid having $\delta^+(V' \setminus \{t\})$ as a minimum
$st$-cut.
\eit
We find a minimum $st$-cut, and we repeat for each $\bar v \in V$. 
Assume that $\{s\} \cup T$ defines one of these minimum cuts. Let $T'=V\cap T$,
$\bar T'=V \setminus T'$. The set $\bar T'$
minimizes $\beta |W| - |E[W]|$ with the constraint that $\bar v \in W$. 
At each iteration of Newton's method the value of $\beta$ increases. Then for
a particular node $\bar v$, we have to solve a sequence of minimum cut problems,
where the capacities of the arcs $(s,v)$ are increasing. Then the parametric network
flow algorithm of \cite{GGT} solves this sequence with the same asymptotic complexity
as one application of the push-preflow algorithm of \cite{GT}. We summarize
below.

\begin{theo}
The arboricity of a hypergraph can be computed with the same asymptotic 
complexity 
of $|V|$ applications of the push-preflow algorithm in a graph with
$O(|V|+|E|)$ nodes. 
\end{theo}

\section{Network Reinforcement}\label{reinf}
The following problem was studied in \cite{Cu}:
Given a graph, a number $k$ and a set of candidate edges, each
of them with an associated cost, find a minimum cost set of candidate edges to be added to the network so it has strength equal to $k$.
Algorithms for this have been given in \cite{Cu}, \cite{GaPar} and \cite{BaReinf}.
Here we modify the algorithm of \cite{BaReinf}
to solve a similar reinforcement problem
in hypergraphs. The algorithm of Section~\ref{separDom} will
be used as a subroutine, so we start with some properties
of its solutions.

 \subsection{Properties of optimal partitions}
The following lemmas are extensions of similar ones for the graphic case, cf.  \cite{BaReinf}.

\begin{lem}\label{l1}
Let $\Phi=\{S_1,\ldots,S_p\}$ be a solution of (\ref{fundSep}), 
and let $\{T_1,\ldots,T_q\}$ be a partition of $S_i$, for
some $i$, $1 \le i \le p$.
Then 
\beqn
  \bar x \big(\delta (T_1,\ldots,T_q)\big) - \beta(q-1) \ge 0.
\eeqn
\end{lem}
\proof If $\bar x \big(\delta(T_1,\ldots,T_q)\big) - \beta(q-1) < 0$
one could improve the solution of (\ref{fundSep}) by 
removing $S_i$ from $\Phi$ and adding $\{T_1,\ldots,T_q\}$.
\endproof

\begin{lem}
\label{l2}
Let $\Phi=\{S_1,\ldots,S_p\}$ be a solution of \eqref{fundSep}, and let
$\{S_{i_1},\ldots,S_{i_l}\}$ be a sub-family of $\Phi$.
Then 
\beqn
  \bar x \big(\delta(S_{i_1},\ldots,S_{i_l})\big) - \beta(l-1) \le 0.
\eeqn
\end{lem}

\proof  If $\bar x \big(\delta(S_{i_1},\ldots,S_{i_l})\big) - \beta (l-1) > 0$,
one could improve the solution of (\ref{fundSep}) by 
removing $\{S_{i_1},\ldots,S_{i_l}\}$ from $\Phi$ and adding their union.
\endproof


\begin{lem}
\label{l3}
Let $\Phi=\{S_1,\ldots,S_p\}$ be a solution of (\ref{fundSep}) in $H$. 
Let $H'$ be the hypergraph obtained by adding one new hyperedge $e$ to $H$.
The following are possible.
\bit
\item[(i)]
If there is an index $i$ such that $e \subseteq S_i$ then
$\Phi$ is a solution of (\ref{fundSep}) in $H'$,
\item[(ii)]
otherwise there exists
a solution $\Phi'$ of (\ref{fundSep}) in $H'$ 
that is either of the form 
%

a) 
$  \Phi'= (\Phi \setminus \{S_i \, : \, i \in I\} )
  \cup \{U=\cup_{i\in I} S_i \},$
for some index set $I \subseteq \{1,\ldots,p\}$, and
$e\in \delta(S_{i_1}, S_{i_2}, \ldots, S_{i_r})$, with $\{i_1, i_2, \ldots, i_r \} \subseteq I$,
see Figure~\ref{fig1}, or

b) $\Phi'=\Phi$. 
\eit
\end{lem}

\proof
(i) Suppose that $e \subseteq S_i$ for some index $i$. If
$\Phi$ is not a solution of \eqref{fundSep} in $H'$, there
is a partition $\Phi'$ such that
$$\bar x(\delta_{H'}(\Phi')) - \beta(|\Phi'|-1)< 
\bar x(\delta_{H'}(\Phi)) - \beta(|\Phi|-1)=\bar x(\delta_{H}(\Phi))-\beta(|\Phi|-1).$$
If $e \notin \delta_{H'}(\Phi)$ we obtain 
$\bar x(\delta_{H}(\Phi')) - \beta(|\Phi'|-1)<\bar x(\delta_{H}(\Phi))-\beta(|\Phi|-1)$,
which is impossible. If $e \in \delta_{H'}(\Phi)$, we have 
$\bar x(\delta_{H}(\Phi')) - \beta(|\Phi'|-1) \le 
\bar x(\delta_{H'}(\Phi')) - \beta(|\Phi'|-1) 
<\bar x(\delta_{H}(\Phi))-\beta(|\Phi|-1)$, a contradiction.

To prove (ii).
Let $\{ T_1,\ldots, T_q \}$ be a solution of (\ref{fundSep}) in $H'$.
Assume that there is a set $S_i$ such that 
$S_i \subseteq \cup_{l=1}^{l=k} T_{j_l}$, $k \ge 2$, and 
$S_i \cap T_{j_l} \ne \emptyset$ for $1 \le l \le k$.
Lemma \ref{l1} implies that 
\beqn
  \bar x \big(\delta_{H}(T_{j_1}\cap S_i,\ldots,T_{j_k}\cap S_i)\big) - \beta (k-1) \ge 0,
\eeqn
and $\bar x \big(\delta_{H'}(T_{j_1},\ldots,T_{j_k})\big) - \beta (k-1) \ge 0$.
Therefore $\{T_{j_1},\ldots,T_{j_k}\}$ can be replaced by their union. So we
can assume that for all $i$ there is an index
$j(i)$ such that $S_i \subseteq T_{j(i)}$.

Now suppose that for some index $j$, $T_j= \cup_{r=1}^{r=l} S_{i_r}$, $l > 1$.
If $e \notin \delta_{H'}(S_{i_1},\ldots,S_{i_l})$,
from Lemma \ref{l2} we have that
\beqn
\bar x \big(\delta_{H'}(S_{i_1},\ldots,S_{i_l})\big) - \beta (l-1)=
  \bar x \big(\delta_{H}(S_{i_1},\ldots,S_{i_l})\big) - \beta (l-1) \le 0,
\eeqn
and we could replace $T_j$ by $\{S_{i_1},\ldots,S_{i_l}\}$. 

If
$e \in \delta_{H'}(S_{i_1},\ldots,S_{i_l})$ and
\beqn
  \bar x \big(\delta_{H'}(S_{i_1},\ldots,S_{i_l})\big) - \beta (l-1) > 0,
\eeqn
we should keep $T_j \in \Phi'$, otherwise we can replace $T_j$ by $\{S_{i_1},\ldots,S_{i_l}\}$. 
\endproof

\begin{figure}[h]
\includegraphics[height=3.5cm, width=8cm]{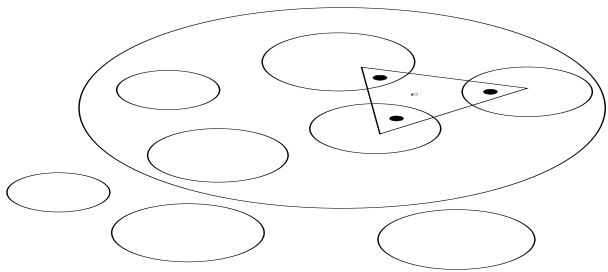} 
\caption{The family $\Phi'$ is obtained by combining some
sets in  $\Phi$.}\label{fig1}
\end{figure}

\subsection{Reinforcement}
We study a slightly different formulation.
For a hypergraph $H=(V,E)$,
 we assume that each hyperedge $e$ has a non-negative per-unit cost $d(e)$ and a non-negative integer capacity $u(e)$, that gives the maximum number of copies allowed of $e$. For a
non-negative integer $k$, we solve the linear program
\beqa
&& \min dx \label{p1}\\
&& x(\delta(\PP)) \ge k (|\PP|-1), \mbox{ for all partitions $\PP$ of $V$,} \label{p2}\\
&& 0 \le x(e) \le u(e). \label{p3}
\eeqa

Its dual is
\beqa
&&\max \sum_\PP \gamma_{\PP} k(|\PP| -1) - \sum u(e) \beta(e)\label{du1}\\
&& \sum_{\PP \, : \, e \in \delta(\PP)} \gamma_\PP \le d(e) + \beta(e), 
\quad \mbox{for all $e$,}\label{du2}\\
&& \gamma \ge 0, \quad \beta \ge 0.\label{du3}
\eeqa
 
 Here for each partition $\PP$ we have a dual variable $\gamma_\PP$.
 Also for each hyperedge $e$ there is a dual variable $\beta(e)$.
 We use a dual algorithm, i.e., constraints \eqref{du2}-\eqref{du3} will always
 be satisfied, and we are going to maximize \eqref{du1}. 
 For the primal problem, constraints \eqref{p3} will
 always be satisfied and \eqref{p2} will be satisfied at the end.
 Complementary slackness will be satisfied at the end.
 We start with an informal description of the algorithm.

 At the beginning we set to zero all dual variables. At each iteration we choose a partition 
$\PP$ and increase the value of $\gamma_\PP$
by $\epsilon$. We have to make sure that constraints \eqref{du2} are not
violated for hyperedges in $\delta(\PP)$. We say that a hyperedge $e$ is {\it tight}
if its constraint \eqref{du2} is satisfied as equation. For a tight hyperedge
$e \in \delta(\PP)$ we have to increase the value of $\beta(e)$ by $\epsilon$.
Let $H'$ be the hypergraph defined by the tight hyperedges, then the dual objective
increases by
\beqn
\epsilon\big(k(|\PP|-1) - u(\delta_{H'}(\PP))\big).
\eeqn
Thus we need to find a partition $\PP$ of $V$ so that $k(|\PP|-1) - u(\delta_{H'}(\PP))$
is positive. Thus
we solve
\beq
\label{part}
\min u(\delta_{H'}(\PP)) - k(|\PP|-1),
\eeq
among all partitions of $V$,
 as in Section~\ref{separDom}.
 Let $\bar \PP=\{S_1,\ldots,S_p\}$ be the partition obtained. 
 Let $(\bar \gamma, \bar \beta)$
 be the current dual solution.
 If the minimum in \eqref{part} is negative
 we look for the largest value of $\epsilon$ so that a new hyperedge becomes tight,
 this is
 \beq\label{eps}
 \bar \epsilon=\min \big\{ d(e) - \sum_{\PP \, : \, e \in \delta(\PP)} 
 \bar \gamma_\PP 
 \, | \,
 e \in \delta_H(\bar \PP) \setminus \delta_{H'}(\bar \PP) \big\}.
 \eeq
 If this minimum is taken over the empty set we say that $\bar \epsilon=\infty$.
 In this case the dual problem is unbounded and the primal problem is
 infeasible.
 
 Let $\bar e$ be an hyperedge giving the minimum in \eqref{eps}. If there is more than
 one hyperedge giving this minimum we pick arbitrarily one to be added to $H'$.
 Then in the next iteration $H'$ will contain a subset of the tight hyperedges, and
 the next value for $\bar \epsilon$ will be zero.
 
 Let $\PP'$ be the solution of \eqref{part} after adding $\bar e$ to $H'$. 
 If $\PP'=\bar\PP$ then $\beta(\bar e)$ could increase and $x(\bar e)$ takes the
 value $u(\bar e)$ to satisfy complementary slackness. We call this {\it Case 1}.
 If $\PP'\ne \bar\PP$ then $\beta(\bar e)$ remains equal to zero and $x(\bar e)$ can 
 take a value less than $u(\bar e)$. This is called {\it Case 2}. The algorithms stops
 when the minimum in \eqref{part} is zero.

 Initially we set $\bar x=0$.
 We have to discuss how to update $\bar x$ in Cases 1 and 2 above.
 
 In Case 1, we set $\bar x(\bar e)=u(\bar e)$. 
 
 In Case 2, from Lemma~\ref{l3} we have that
 \beqn
 \PP'=(\bar \PP \setminus \{S_i \, : \, i \in I \}) \cup \{U=\cup_{i \in I} S_i\},
 \eeqn
for some index set $I \subseteq \{1, \ldots, p\}$, and 
$\bar e \in \delta(S_{i_1},\ldots,S_{i_r})$, with $\{i_1,\ldots,i_r\} \subseteq I$.
Let $\PP_I=\{ S_i \, | \, i \in I\}$, 
Lemma~\ref{l1}
implies
\beq \label{eqxxx}
u(\delta_{H'}(\PP_I) ) - k(|I|-1) \ge 0.
\eeq

At this point we have
$\bar x(e)=u(e)$ for $e \in \delta_{H'}(\PP_I) \setminus \{\bar e\}$, and
$\bar x(e)=0$ for all other hyperedges $e \in \delta_{H}(\PP_I)$.
If $0 < \bar x(e) < u(e)$ for some edge $e$, then the value  $\bar x(e)$ has been
set in Case 2 in some former iteration, and $e \subseteq S_i$ for some set $S_i \in \bar \PP$.
Let 
 \beqn
 \lambda=k(|I|-1) - \bar x( \delta_H(\PP_I)\setminus \{\bar e\}),
 \eeqn
 Inequality \eqref{eqxxx} implies
 $\lambda \le u(\bar e)$.
Then we set  
$\bar x(e)=\lambda$, and we have $\bar x( \delta_H(\PP_I))=k(|I|-1)$.

\begin{lem}\label{case2}
For the set $U$ defined in Case 2, we have 
\bit
\item[(a)] $\bar x(E[U])=k (|U|-1)$, and
\item[(b)] $\bar x(\delta_H(T_1, \ldots, T_q)) \ge k(q-1)$ for each partition $\{T_1, \ldots, T_q\}$
of $U$.
\eit
\end{lem} 
\proof
We use induction, so we assume that the lemma holds for each
$S_i$, $i \in I$.

Consider (a).
We have $\bar x(E[S_i])=k(|S_i|-1)$ for $i \in I$, by the induction
hypothesis. And 
since $\bar x(\delta_H(\PP_I))=k(|I|-1)$, we have
$\bar x(E[U])=k (|U|-1)$.

To prove (b) assume that $\{U_1,\ldots, U_r\}$ is a solution
of 
\beq\label{partY}
\min \bar x(\delta_H(\mathscr{Q})) - k (|\mathscr{Q})|-1)
\eeq
over all partitions $\mathscr{Q}$ of $U$. 
The induction hypothesis implies that for each $i \in I$,
$S_i \subseteq U_{j(i)}$ for
some index $j(i)$.
If $\bar x (\delta_H(U_1, \ldots, U_r) ) - k(r-1) < 0$,
and since 
\hbox{$\bar x(\delta_H(\PP_I))-k(|I|-1)=0$},
there is a family $\{S_{i_1}, \ldots, S_{i_m}\} \subset \PP_I$, such that
$\cup_{j=1}^{j=m} S_{i_j} = U_l$ for some index $l$,
and 
\beqn
\bar x(\delta_H(S_{i_1}, \ldots, S_{i_m})) > k(m-1).
\eeqn
This contradicts Lemma~\ref{l2} and the definition of $\PP'$.
\endproof

The first partition of $V$ consists of all singletons.
At each iteration a new hyperedge is added to $H'$.
 In some cases some sets in the
 family $\bar \PP$ are combined into one. Each time that several sets
 of a partition are combined into one set $U$, the update of $\bar x$
 implies $\bar x(E[U])=k(|U|-1)$, as shown in Lemma~\ref{case2}.
 When the minimum in \eqref{part} is zero, then $\PP=\{V\}$ is
 a solution. Thus
 $U=V$,
 and we have $\bar x(E)=k(|V|-1)$. Then Lemma~\ref{case2} shows that
 we have a primal feasible solution. 
 
 At the end, consider a partition $\PP=\{S_1, \ldots, S_p\}$
 with $\bar \gamma_\PP > 0$.
 We have $\sum_i \bar x(E[S_i])+ \bar x(\delta_H(\PP))=\bar x(E)=k(|V|-1)$.
Lemma~\ref{case2} shows that $\bar x(E[S_i]) = k(|S_i|-1)$, for all $i$.
Therefore $\bar x(\delta_H(\PP))=k(|V|-1) -\sum_i k(|S_i|-1)=
k(|\PP|-1)$.
 Thus at the end the vectors
 $\bar x$ and $(\bar \gamma, \bar \beta)$ satisfy the complementary
 slackness conditions. Now we give a formal description of the algorithm.

 \vskip 0.5 cm
\centerline{\bf Reinforcement}
\begin{itemize}
\item{{\bf Step 0}.} 
Start with $\bar \gamma=0$, $\bar \beta=0$, 
$\bar x=0$,
$\bar d(e)= d(e)$ for all $e \in E$.
$\bar \PP$ consisting of all singletons, and $H'=(V, \emptyset)$.
\item{{\bf Step 1}.} Compute 
\beq \label{epsBar}
\bar \epsilon = \min \{ \bar d(e) \, | \, e \in 
\delta_H(\bar \PP) \setminus \delta_{H'}(\bar \PP) \}.
\eeq

If $\bar \epsilon=\infty$ stop, the problem is
infeasible. \newline
Otherwise update 
$\bar \beta(e) \leftarrow 
\bar \beta(e) +  \bar \epsilon$ for all $e \in \delta_{H'}(\bar \PP)$,
\newline
$\bar \gamma_{\bar \PP} \leftarrow \bar \gamma_{\bar \PP} + \bar \epsilon$,
\newline
$\bar d(e) \leftarrow 
\bar d(e) - \bar \epsilon$ for all 
$e \in \delta_H(\bar \PP) \setminus \delta_{H'}(\bar \PP)$.
\item{{\bf Step 2}.} Let $\bar e$ be a hyperedge giving the minimum in
\eqref{epsBar}, add $\bar e$ to $H'$. Solve problem \eqref{part} in $H'$
to obtain a partition $\PP'$.
\item{{\bf Step 3}.} If $\bar \PP=\PP'$ update $\bar x$ as
in Case 1. Otherwise update as in Case 2. 
If the minimum in \eqref{part} is zero 
stop, 
otherwise set $\bar \PP \leftarrow \PP'$ and go to Step 1.
\end{itemize}
 
 This algorithm takes at most $|E|$ iterations, where each of them requires
to solve \eqref{part}. 
For the case when $H$ is a graph,
$k=1$, and $u(e)=\infty$ for every edge $e$, this algorithm is similar
to Kruskal's algorithm for minimum spanning trees \cite{krus}.
Now we summarize the results of this section.

\begin{theo}
The reinforcement problem reduces to $|E||V|$ minimum cut problems
in a graph with $O(|V|+|E|)$ nodes.
\end{theo}

 \begin{theo}
 The vector $\bar x $ is an optimal solution of \eqref{p1}-\eqref{p3}. If $k$ is a nonnegative
integer and the capacities $u$ are integer, then $\bar x$  is integer valued. 
Also if $d$ is integer
valued then there is an optimal solution of \eqref{du1}-\eqref{du3} that is also integer valued. Thus
the system \eqref{p2}-\eqref{p3} is totally dual integral.
 \end{theo}
 \endproof

\section{Final Remarks}
We have given extensions to hypergraphic matroids of several algorithms that
had been developed for graphs.  This work is the basis of an algorithm for packing
hypertrees, and its uses to compute lower and upper bounds for the $k$-cut
problem in hypergraphs. This will be presented in a forthcoming paper
\cite{packHyper}.

\bibliographystyle{siam}
\bibliography{hyperg.bib}



\end{document}